\documentclass[11pt]{amsart}
\usepackage{txfonts}
\usepackage{pifont}
\usepackage{amsfonts}
\usepackage{amssymb}
\usepackage{mathrsfs}
\usepackage[abs]{overpic}
\usepackage{cases}
\usepackage{amsfonts}
\usepackage{amsmath}
\usepackage{mathrsfs}
\usepackage{amssymb}
\usepackage{amsmath,color}
\usepackage{amsthm}
\usepackage{amstext}
\usepackage{amsopn}
\usepackage{texdraw}
\usepackage{graphicx}
\usepackage{extarrows}
\usepackage{amsmath}
\usepackage{amsthm}
\usepackage{amstext}
\usepackage{amsopn}
\usepackage{graphicx}
\usepackage{subfigure}
\usepackage{cite}
\usepackage{bbm}
\newtheorem{theorem}{Theorem}[section]
\newtheorem{lemma}[theorem]{Lemma}
\newtheorem{corollary}[theorem]{Corollary}
\newtheorem{proposition}[theorem]{Proposition}

\newtheorem{problem}[theorem]{Problem}

\newtheorem{claim}[theorem]{Claim}

\theoremstyle{definition}

\newtheorem{example}[theorem]{Example}

\theoremstyle{remark}
\newtheorem{remark}[theorem]{Remark}

\newcommand{\eqb}{\begin{equation}}
\newcommand{\eqe}{\end{equation}}
\newcommand{\spb}{\begin{split}}
\newcommand{\spe}{\end{split}}
\newcommand{\cab}{\begin{cases}}
\newcommand{\cae}{\end{cases}}

\newcommand{\thmb}{\begin{theorem}}
\newcommand{\thme}{\end{theorem}}
\newcommand{\prob}{\begin{proposition}}
\newcommand{\proe}{\end{proposition}}
\newcommand{\leb}{\begin{lemma}}
\newcommand{\lee}{\end{lemma}}
\newcommand{\cob}{\begin{corollary}}
\newcommand{\coe}{\end{corollary}}
\newcommand{\enb}{\begin{enumerate}}
\newcommand{\ene}{\end{enumerate}}
\newcommand{\mb}{\begin{matrix}}
\newcommand{\me}{\end{matrix}}
\newcommand{\cenb}{\begin{center}}
\newcommand{\cene}{\end{center}}

\newcommand{\prb}{\begin{proof}}
\newcommand{\pre}{\end{proof}}
\newcommand{\rb}{\begin{remark}}
\newcommand{\re}{\end{remark}}
\newcommand{\lf}{\lfloor}
\newcommand{\rf}{\rfloor}

\newcommand{\lt}{\left}
\newcommand{\rt}{\right}
\newcommand{\la}{\langle}
\newcommand{\ra}{\rangle}
\makeatletter
\newcommand{\vast}{\bBigg@{4.5}}
\makeatother

\numberwithin{equation}{section}

\begin{document}

 \title[Circularly invariant uniformizable probability measures]
  {Circularly invariant uniformizable probability measures for linear transformations}

\author{Chuang Xu}

\address{
Department of Mathematical and Statistical Sciences,
University of Alberta, Edmonton, Alberta,
T6G 2G1 Canada.}

\email{cx1@ualberta.ca.}




\noindent

\begin{abstract}
In this paper, we prove a threshold result on the existence of a circularly invariant uniformizable probability measure (CIUPM) for linear transformations with non-zero slope on the line. We show that there is a threshold constant $c$ depending only on the slope of the linear transformation such that there exists a CIUPM if and only if its support has a diameter at least as large as $c.$ Moreover, the CIUPM is unique up to translation if the diameter of the support equals $c.$
\end{abstract}  \maketitle
\section{Introduction}

This paper investigates a variant of invariant probability measure for linear transformations on the line. Let $\mathbb{R}$ be the set of real numbers and $\mathbb{T}$ the unit circle identified with $[0,1[$ via the canonical mapping $t\mapsto e^{2\pi{\rm i}t},$ where the usual notations for intervals, e.g., $[a,b[\ =\{x\in\mathbb{R}: a\le x<b\}$  are used throughout. Let $\lambda$ and $\lambda_{\mathbb{T}}$ be the Lebesgue measure on $\mathbb{R}$ and $\mathbb{T},$ respectively. Denote by $R:\cab\mathbb{R}\to\mathbb{T}\\ x\mapsto\la x\ra\cae$ the rotation mapping with $\la x\ra$ being the fractional part of $x,$ and ${\rm diam}(A):=\underset{x,y\in A}{\sup}|x-y|$ the diameter of a set $A\subset\mathbb{R}.$ Note that ${\rm diam}(A)\ge\lambda(A)$ for every set $A$ and the equality holds if and only if $\lambda(A)=+\infty$ or $A$ is an interval but a Lebesgue measure zero set (i.e., $\lambda([\inf A,\sup A]\setminus A)=0$ with $\inf A$ and $\sup A$ denoting the infimum and supremum of $A,$ respectively). For every continuous monotone transformation $T:\mathbb{R}\to\mathbb{R},$ let $\mu\circ T^{-1}$ be the {\em induced} (or {\em push-forward}) probability measure for $T,$ and $\la \mu\ra$ denotes $\mu\circ R^{-1}$ for convenience. Note that $\la\mu\ra$ and $\lt\la\mu\circ T^{-1}\rt\ra$ both are probability measures on $\mathbb{T}.$ A measure $\mu$ on $\mathbb{R}$ is a {\em a circularly invariant uniformizable probability measure} (CIUPM) if $$\la\mu\ra=\lt\la\mu\circ T^{-1}\rt\ra=\lambda_{\mathbb{T}}.$$ Obviously, every CIUPM is absolutely continuous (w.r.t. $\lambda.$)

Our motivation for the study of CIUPMs comes from {\em uniform distribution theory}. For every sequence $(x_n)$ of real numbers, let \[\mu_n:=\mu_n(x_n)=\frac{1}{n}\sum_{i=1}^n\delta_{x_i}\]be the sequence of probability measures on $\mathbb{R}$ generated by $(x_n).$ It is known that for some convex monotone transformations $T:\mathbb{R}\to\mathbb{R}$ like the exponential transformation or the trivially convex linear transformation, there exists a {\em uniformly distributed modulo one} (u.d. {\rm mod} 1) sequence $(x_n)$ of real numbers such that $(T(x_n))$ is also u.d. {\rm mod} 1. (For the definition of u.d. {\rm mod} 1 sequences, cf. \cite{KN}.) Precisely, for some convex monotone transformation $T:\mathbb{R}\to\mathbb{R},$ there exists a sequence $(x_n)$ such that $\lt\la\mu_n(x_n)\rt\ra\to\lambda_{\mathbb{T}},\ \lt\la\mu_n(T(x_n))\rt\ra\to\lambda_{\mathbb{T}}$
weakly. In fact, it follows directly from Weyl's criterion \cite[Chapter~1, Theorem~2.1]{KN} that for every linear transformation with a non-zero slope, $(x_n)$ is u.d. {\rm mod} 1 if and only if $(T(x_n))$ is u.d. {\rm mod} 1. Also, for $T(x)=e^x,$ both $(\alpha n)$ and $(T(\alpha n))$ are u.d. {\rm mod} 1, for almost all irrational numbers $\alpha$ \cite[Chapter~1, Corollary~4.1]{KN}. However, it remains open whether $\lt(a^n\rt)$ is u.d. {\rm mod} 1, for some specific positive number $a,$ for instance, when $a=e,\ \pi$ or even as simple as $3/2$ \cite[p.36]{KN}.

Then a natural analogous question arises: for a given convex monotone transformation $T:\mathbb{R}\to\mathbb{R},$ does there always exist a probability measure $\mu$ on $\mathbb{R}$ such that $\la\mu\ra=\lt\la\mu\circ T^{-1}\rt\ra=\lambda_{\mathbb{T}}\ \text{?}$
In other words, does there exist a CIUPM for $T$ ? As we will show in this paper, though in the discrete version it is trivial that for any u.d. {\rm mod} 1 sequence $(x_n),$ $(T(x_n))$ is u.d. {\rm mod} 1 for any linear $T$ (with a non-zero slope), it may not be as trivial to show the existence of a CIUPM for a linear transformation $T$ as shown later in Section~3.

This work, as a first try, answers the question for (the trivially convex) linear transformations $T.$ For nonlinear convex transformations, like the exponential transformation, the problem is more difficult in that such CIUPMs, if they exist, cannot be easily solved by their densities as for the linear case in this paper. Indeed, even for a piecewise linear transformation (for instance $T(x)=x\mathbbm{1}_{]-\infty,0]}+\sqrt{2}x\mathbbm{1}_{]0,+\infty[}$), situation becomes much more complicated than the linear case. This will be illustrated more clearly when solving the equations for the densities of a CIUPM in the proof of the main result. Except for the existence of a CIUPM, we present a threshold result characterizing how \lq\lq slim\rq\rq\ a CIUPM can be: For a linear transformation $T$ with non-zero slope, there exists a CIUPM $\mu$ for $T$ if and only if ${\rm diam}({\rm supp}\ \mu)\ge c$ for some positive constant $c$ depending only on the slope of the linear transformation, where ${\rm supp}\ \mu$ is the support of $\mu$ (i.e., the smallest closed subset in $\mathbb{R}$ of full $\mu$ measure). Moreover, the CIUPM is unique up to translation if ${\rm diam}\lt({\rm supp}\ \mu\rt)=c.$

Let us mention some related works on invariant measures for \lq\lq almost\rq\rq\ linear transformations on $[0,1].$ Kopf \cite{Ko} gave a formula for the densities of invariant measures for piecewise linear transformations on $[0,1].$ G\'{o}ra \cite{G1,G2} found an explicit formula for the densities of  invariant measures for arbitrary eventually expanding piecewise linear transformations whose slopes are not necessarily the same on $[0,1].$

For $\alpha\in\mathbb{R},$ $\beta\neq0,$ define $\lt\la T_{\alpha,\beta}\rt\ra:\cab\mathbb{T}\to\mathbb{T},\\ x\mapsto\lt\la\beta x+\alpha\rt\ra.\cae$ For $0\le\alpha<1$ and $\beta>1,$ Parry \cite{P2} gave an explicit formula for the unique invariant measure. Halfin \cite{Ha} showed this invariant measure is positive. Hofbauer \cite{H0,H2,H4,H6} proved that this measure is absolutely continuous (w.r.t. $\lambda$), its entropy equals $\log\beta,$ and its support is a finite union of intervals; he also showed the uniqueness of invariant measures with maximal entropy and determined the region of $(\beta,\alpha)$-plane where ${\rm supp}\ \mu\subset[0,1].$ Faller and Pfister \cite{FP} studied normal points for $\lt\la T_{\alpha,\beta}\rt\ra.$

\section{Preliminaries}

The following standard notations are used throughout. The integers, natural and rational numbers are denoted by $\mathbb{Z},$ $\mathbb{N}$ and $\mathbb{Q},$ respectively. For any real number $x,$ denote by $\lf x\rf$ the floor of $x,$ i.e., the largest integer not exceeding $x,$ and hence $x=\lf x\rf+\la x\ra.$ For $x\in\mathbb{R}$ and $A\subset\mathbb{R},$ let $A+x=\{y+x:\ y\in A\}.$ If $A\subset\mathbb{R}$ is an interval and $f: A\to\mathbb{R}$ is monotone, then  $f(a-)=\lim_{\varepsilon\downarrow0}f(a-\varepsilon)$ and $f(a+)=\lim_{\varepsilon\downarrow0}f(a+\varepsilon)$ both exist for every interior point $a$ of $A.$

Recall that two integers $p,\ q$ are {\em coprime} if they have one as their {\em greatest common divisor} \cite[p.5]{H}.  For every $\beta\in\mathbb{Q}\setminus\{0\},$ let $(p_\beta,q_{\beta})$ be the unique pair of coprime positive integers such that $|\beta|=p_{\beta}/q_{\beta}$ and let $s_{\beta}:=q_{\beta}\la|\beta|\ra.$ Note that for $\beta\notin\mathbb{N},$ $1\le s_{\beta}\le q_{\beta}-1$ is an integer coprime with $q_{\beta}.$

For a complete metric space $X$($=\mathbb{R}$ or $\mathbb{T}$), let $\mathcal{P}(X)$ be the family of all Borel probability measures on $X.$
For $\nu\in\mathcal{P}(\mathbb{T}),$ define its associated {\em distribution function} as\[F_{\nu}(t)=\nu([0,t]),\ t\in\mathbb{T}.\] Let $F_{\la\mu\ra},$ $F_{\mu}$ be the distribution functions and, $\rho_{\la\mu\ra},$ $\rho_{\mu}$ the densities (if they exist) of $\la\mu\ra$ and $\mu,$ respectively. For $c\ge0,$ let $S_c=\lt\{\mu\in\mathcal{P}(\mathbb{R}):{\rm diam}({\rm supp}\ \mu)= c\rt\}.$ For a (piecewise) continuous monotone transformation $T$ on $\mathbb{R},$ let $U_T$ be the set of all CIUPMs for $T$ in $\mathcal{P}(\mathbb{R}).$ 

For the transformation $T_{\alpha,\beta}:\mathbb{R}\to\mathbb{R}$ defined by \eqb\label{0-0}T_{\alpha,\beta}(x)=\beta x+\alpha,\eqe
we study in the next section the problem below on the existence and shortest ``length'' of a CIUPM: \begin{problem}\label{p0} For $\alpha\in\mathbb{R}$ and $\beta\neq0,$ what is $\inf\ \lt\{c:\ U_{T_{\alpha,\beta}}\cap S_c\neq\emptyset\rt\}$ ?
\end{problem}
In the next section, we prove a threshold result on the existence of CIUPMs: For any $\alpha\in\mathbb{R},$ $\beta\neq0,$ there exists $c_{\beta}>0$ such that $U_{T_{\alpha,\beta}}\cap S_c\neq\emptyset$ if and only if $c\ge c_{\beta}.$ Moreover, $U_{T_{\alpha,\beta}}\cap\ S_{c_{\beta}}=\lt\{\mu_{\beta}\circ T_{\gamma,1}^{-1}\rt\}_{\gamma\in\mathbb{R}}$ for some $\mu_{\beta}\in\mathcal{P}(\mathbb{R});$ for every $c>c_{\beta},$ there exist $\mu_1,\mu_2 \in U_{T_{\alpha,\beta}}\cap\ S_{c}$ such that $\mu_1\neq\mu_2\circ T_{\gamma,1}^{-1},\ \forall\ \gamma\in\mathbb{R}.$ In other words, by Proposition~\ref{le3} (i) below, such CIUPM is unique up to translation in $S_{c_{\beta}}$ while {\em not} unique in $S_c$ for $c>c_{\beta}.$

Now we give some preliminary results for the proof of the unique main result in the next section.

By definition, it follows from the Radon-Nikodym theorem (cf. \cite[p.158]{K}) that:
\prob\label{le4}
Assume $\mu\in\mathcal{P}(\mathbb{R}).$ If $\la\mu\ra=\lambda_{\mathbb{T}},$ then $\mu$ is absolutely continuous (w.r.t. $\lambda$), and thus its density $\rho_{\mu}$ exists and is finite ($\lambda$-)almost everywhere (a.e.). In particular, any CIUPM is absolutely continuous with a.e. finite density.
\proe
The result below follows from the definition of the Perron-Frobenius operator (cf. \cite[p.42]{LM}).
\prob\label{le2}
For every $\mu\in\mathcal{P}(\mathbb{R}),$
\[
    F_{\la\mu\ra}(t)=\sum_{k\in\mathbb{Z}}\Big[F_{\mu}(t+k)-F_{\mu}(k-)\Big], \forall\ t\in\mathbb{T}.
\] If $\rho_{\mu}$ exists a.e., then $\rho_{\la\mu\ra}$ exists a.e. on $\mathbb{T},$ moreover,
\[
    \rho_{\la\mu\ra}(t)=\sum_{k\in\mathbb{Z}}\rho_{\mu}(t+k),\ \text{a.e.\ on}\ \mathbb{T}
\]
\proe
\rb
By straightforward calculations for densities via Proposition~\ref{le2}, it is not difficult to see that for any $\mu\in\mathcal{P}(\mathbb{R})$ such that $\la\mu\ra=\lt\la\mu\circ T_{\alpha,\beta}^{-1}\rt\ra,$ it is not necessarily that $\mu=\mu\circ\lt\la T_{\alpha,\beta}\rt\ra^{-1}$ except for the trivial cases of $\alpha=0$ and $\beta\in\mathbb{N}.$
\re
We leave the proofs of the following properties of $U_{T_{\alpha,\beta}}$ an exercise to the reader.
\prob\label{le3}
Assume $\beta\neq0.$\\
\noindent
{\rm (i)} (Translation invariance.) For every $\alpha\in\mathbb{R},$ $U_{T_{0,\beta}}=U_{T_{\alpha,\beta}};$\\
\noindent
{\rm (ii)} (Convexity.) The set $U_{0,\beta}$ is convex: for every $n\in\mathbb{N},$ let $p=(p_1,\cdots,p_n)$ be a probability vector. If $\mu_i\in U_{T_{0,\beta}}$ for all $i=1,\cdots,n,$ then $\sum_{i=1}^np_i\mu_i\in U_{T_{0,\beta}};$\\
\noindent
{\rm (iii)} $\mu\in U_{T_{0,\beta}}$ if and only if $\mu\circ T_{0,-1}^{-1}\in U_{T_{0,-\beta}};$\\
\noindent
{\rm (iv)} $\mu\in U_{T_{0,\beta}}$ if and only if $\mu\circ T_{0,\beta}^{-1}\in U_{T_{0,1/\beta}}.$
\proe

From Proposition~\ref{le3}, it easily follows:
\cob\label{co1}
Assume $\beta\neq0,$ $\alpha\in\mathbb{R}$ and $c\ge0.$ Suppose $U_{T_{\alpha,\beta}}\cap S_c\neq\emptyset.$ Then for every $\widetilde{c}>c,$ $U_{T_{\alpha,\beta}}\cap S_{\widetilde{c}}\neq\emptyset;$ moreover, $$\lt\{\mu\circ T_{\gamma,1}^{-1}\rt\}_{\gamma\in\mathbb{R}}\subsetneqq U_{T_{\alpha,\beta}}\cap S_{\widetilde{c}},\ \forall\ \mu\in U_{T_{\alpha,\beta}}\cap S_{\tilde{c}},$$
\coe
\prb
We prove it by construction. Let $\nu\in U_{T_{\alpha,\beta}}\cap S_c.$ It follows from Proposition~\ref{le3}(i) that $\nu\circ T_{\tilde{c}-c,1}^{-1}\in U_{T_{\alpha,\beta}}.$ Let $\widetilde{\nu}=\frac{1}{2}\nu+ \frac{1}{2}\nu\circ T_{\tilde{c}-c,1}^{-1}.$ Then $\widetilde{\nu}\in U_{T_{\alpha,\beta}}$ by Proposition~\ref{le3}(i) and (ii). Moreover, ${\rm diam}\lt({\rm supp}\ \widetilde{\nu}\rt)={\rm diam}\lt({\rm supp}\ \nu\rt)+\tilde{c}-c=\tilde{c}.$ This yields $\widetilde{\nu}\in U_{T_{\alpha,\beta}}\cap S_{\tilde{c}}.$ Analogously, $\widetilde{\mu}=\frac{1}{3}\nu+ \frac{1}{3}\nu\circ T_{(\tilde{c}-c)/2,1}^{-1}+\frac{1}{3}\nu\circ T_{\tilde{c}-c,1}^{-1}\in U_{T_{\alpha,\beta}}\cap S_{\tilde{c}}.$ This yields $\lt\{\mu\circ T_{\gamma,1}^{-1}\rt\}_{\gamma\in\mathbb{R}}\neq U_{T_{\alpha,\beta}}\cap S_{\tilde{c}},$ for all $\mu\in U_{T_{\alpha,\beta}}\cap S_{\tilde{c}}.$ On the other hand, it trivially follows from Proposition~\ref{le3}(i) that for $\lt\{\mu\circ T_{\gamma,1}^{-1}\rt\}_{\gamma\in\mathbb{R}}\subset U_{T_{\alpha,\beta}}\cap S_{\tilde{c}},$ for all $\mu\in U_{T_{\alpha,\beta}}\cap S_{\tilde{c}}.$
\pre
Two real numbers $x$ and $y$ are {\em rationally independent} if one is a rational multiple of the other, i.e., equation $r_1x+r_2y=0$ only admits the trivial solution $r_1=r_2=0$ in $\mathbb{Q}^2.$

The following elementary result on rational independence of real numbers follows directly from a Chebyshev's theorem on \cite[p.266]{H}.

\prob\label{le1}
Assume $x_1,\ x_2\in\mathbb{R}.$ If $x_1,\ x_2$ are rationally independent, then the sequence $\lt(mx_1+nx_2\rt)_{m,n\in\mathbb{Z}}$ is dense in $\mathbb{R}.$
\proe

\section{An Answer to Problem~\ref{p0}}

In this section, we give an answer to Problem~\ref{p0} via a threshold result on the existence of CIUPMs for linear transformations $T_{\alpha,\beta}$ defined by \eqref{0-0}.

\begin{figure}[h]
\centering
    \subfigure[]{\begin{overpic}[width=1.3in]{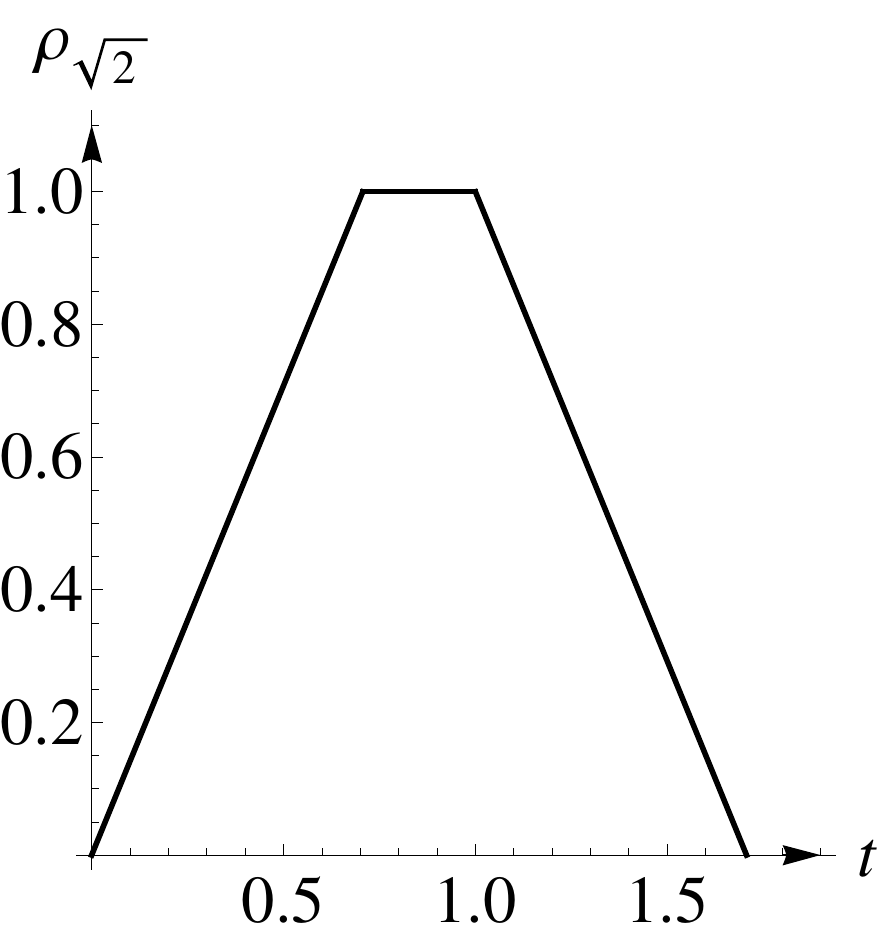}
\put(5,4){\makebox(0,0)[cc]{\tiny$0.0$}}
\end{overpic}} \hspace{0.45cm}
 \subfigure[]{\begin{overpic}[width=1.3in]{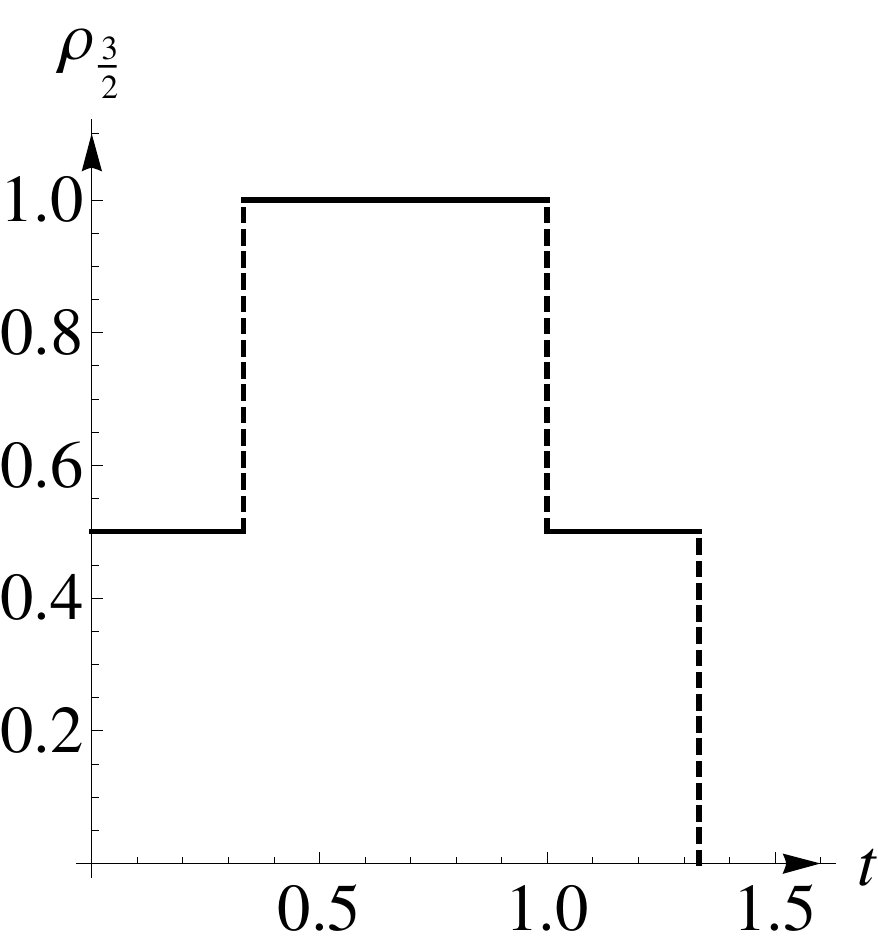}
\put(5,4){\makebox(0,0)[cc]{\tiny$0.0$}}
\end{overpic}}\\ \hspace{0.45cm}
\subfigure[]{\begin{overpic}[width=1.3in]{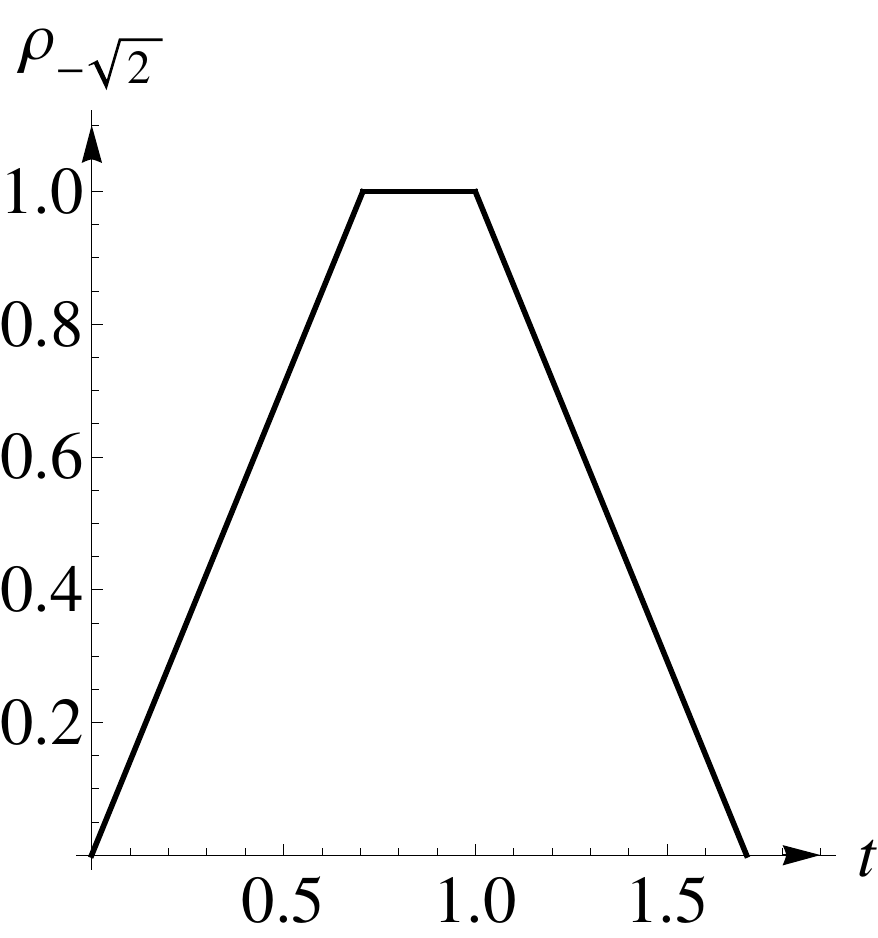}
\put(5,4){\makebox(0,0)[cc]{\tiny$0.0$}}
\end{overpic}} \hspace{0.45cm}
\subfigure[]{\begin{overpic}[width=1.3in]{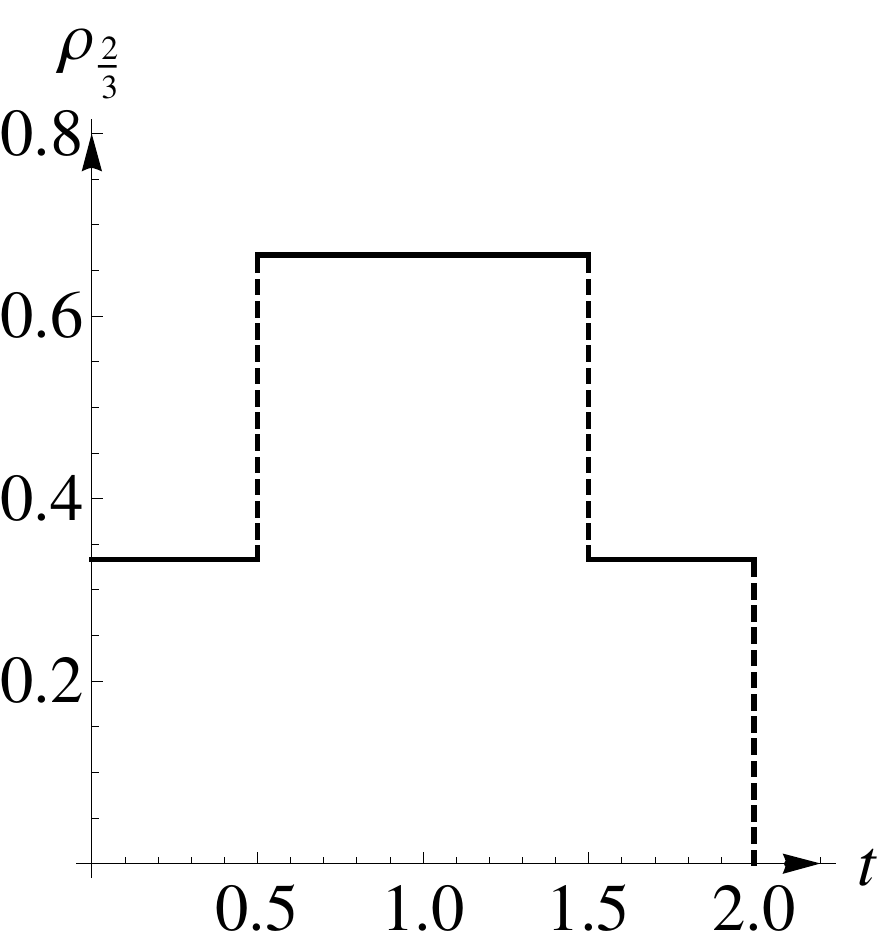}
\put(5,4){\makebox(0,0)[cc]{\tiny$0.0$}}
\end{overpic}}
\caption{\small Profiles for $\rho_{\beta}.$}\label{f0}
\end{figure}

Before stating the threshold result, let us look at two simple examples, which may give some intuitive picture of a \lq\lq slimmest\rq\rq\ CIUPM.
Consider first a linear transformation with an irrational slope.
\begin{example}
Let $T_{0,\sqrt{2}}(x)=\sqrt{2}x.$ By Proposition~\ref{le2}, simple calculations verify that $\mu$ is a CIUPM for $T_{0,\sqrt{2}}$ with density $\rho_{\mu}(t)=\begin{cases}
      \sqrt{2}t\ \ \ \ \ \ \ \ \ \ \ \ \ \ \ \ \ \ \ \text{if}\ t\in\lt[0,1/\sqrt{2}\rt[,\\
      1\ \ \ \ \ \ \ \ \ \ \ \ \ \ \ \ \ \ \ \ \ \ \ \ \text{if}\ t\in\lt[1/\sqrt{2},1\rt[,\\
      -\sqrt{2}t+1+\sqrt{2}\ \text{if}\ t\in\lt[1,1+1/\sqrt{2}\rt],\\
      0\ \ \ \ \ \ \ \ \ \ \ \ \ \ \ \ \ \ \ \ \ \ \ \ \ \text{elsewhere}.
    \end{cases}$ Notice that ${\rm supp}\ \mu$ is an interval and ${\rm diam}\lt({\rm supp}\ \mu\rt)=\lambda\lt({\rm supp}\ \mu\rt)=1+\frac{1}{\sqrt{2}}.$ See also Fig.~\ref{f0} (a).
\end{example}
Now we turn to a linear transformation with a rational slope.\begin{example}
Let $T_{0,3/2}(x)=3x/2.$ It also follows from Proposition~\ref{le2} that $\mu$ is a CIUPM for $T_{0,3/2}$ with density $\rho_{\mu}(t)=\begin{cases}
     \frac{1}{2}\ \text{if}\ t\in\lt[0,1/3\rt[\bigcup\lt[1,4/3\rt[,\\
      1\ \text{if}\ t\in\lt[1/3,1\rt[,\\
       0\ \text{elsewhere}.
    \end{cases}$ Note that ${\rm supp}\ \mu$ is an interval and ${\rm diam}\lt({\rm supp}\ \mu\rt)=\lambda\lt({\rm supp}\ \mu\rt)=4/3<1+\frac{1}{3/2}.$
\end{example}
From the above two examples, one may expect there always exists a CIUPM for every linear transformation with a non-zero slope. In fact, as illustrated by the following main result, these two CIUPMs in these examples are the \lq\lq slimmest\rq\rq\ (in the sense that the support of CIUPM has the shortest diameter).

For every $\beta\neq0,$ define a constant \[ c_{\beta}=\begin{cases}
    1+\frac{1}{|\beta|}-\frac{1}{p_{\beta}}\ \text{if}\ \beta\in\mathbb{Q},\\
    1+\frac{1}{|\beta|}\ \ \ \ \ \ \ \ \ \text{if}\ \beta\in\mathbb{R}\setminus\mathbb{Q},
  \end{cases}\] and a probability measure $\mu_{\beta}$ by its density $\rho_{\beta}:$\\
{\rm(i)}    \eqb\label{0-1}
    \rho_{\beta}(t)=\begin{cases}
      \beta t\ \ \ \ \ \ \ \ \ \ \ \ \ \ \ \ \text{if}\ t\in\lt[0,\frac{1}{\beta}\rt[,\\
      1\ \ \ \ \ \ \ \ \ \ \ \ \ \ \ \ \ \text{if}\ t\in\lt[\frac{1}{\beta},1\rt[,\\
      -\beta t+1+\beta\ \text{if}\ t\in\lt[1,1+\frac{1}{\beta}\rt],\\
      0\ \ \ \ \ \ \ \ \ \ \ \ \ \ \ \ \ \text{elsewhere},
    \end{cases}
  \eqe if $\beta\in[1,+\infty[\setminus\mathbb{Q};$\\
\noindent{\rm(ii)}  \eqb\label{0-3}
    \rho_{\beta}(t)=\begin{cases}
      \frac{j}{q_{\beta}}\ \text{if}\ t\in\lt[\frac{j-1}{p_{\beta}},\frac{j}{p_{\beta}}\rt[\bigcup\lt[1+
      \frac{q_{\beta}-j-1}{p_{\beta}},1+\frac{q_{\beta}-j}{p_{\beta}}\rt[,\\ \ \ \ \ \ \ \text{for}\ j=1,\cdots,q_{\beta}-1,\\
      1\ \ \text{if}\ t\in\lt[\frac{q_{\beta}-1}{p_{\beta}},1\rt[,\\
      0\ \ \text{elsewhere},
    \end{cases}\eqe if $\beta\in[1,+\infty[\ \cap\ \mathbb{Q};$\\
\noindent
{\rm(iii)}  $\rho_{\beta}(\cdot)=\beta\rho_{1/\beta}(\beta \cdot)$ if $\beta\in]0,1[;$ \\
\noindent
{\rm(iv)} $\rho_{\beta}(\cdot)=\rho_{-\beta}(\cdot)$ if $\beta\in]-\infty,0[.$
\begin{theorem}\label{th1}
For every $\alpha\in\mathbb{R},\ \beta\neq0,$ let $T_{\alpha,\beta}$ be defined as in \eqref{0-0}.  Then $U_{T_{\alpha,\beta}}\cap S_c\neq\emptyset$ if and only if $c\ge c_{\beta};$ moreover, $U_{T_{\alpha,\beta}}\cap S_{c_{\beta}}=\lt\{\mu_{\beta}\circ T_{\gamma,1}^{-1}\rt\}_{\gamma\in\mathbb{R}},$ but for every $c>c_{\beta},$ $$U_{T_{\alpha,\beta}}\cap S_{c}\neq\lt\{\mu\circ T_{\gamma,1}^{-1}\rt\}_{\gamma\in\mathbb{R}},\ \forall\ \mu\in U_{T_{\alpha,\beta}}\cap S_{c},$$ where $c_{\beta}$ and $\mu_{\beta}$ are defined as above.
\end{theorem}
\prb
Throughout denote for short $T_{0,\beta}$ and $U_{0,\beta}$ by $T_{\beta}$ and $U_{\beta},$ respectively. Note that the conclusions are trivial for $\beta\in\mathbb{N}.$ By Proposition~\ref{le3}, Corollary~\ref{co1}, and the symmetry $\rho_{\beta}(t)=\rho_{\beta}(c_{\beta}-t),$ it suffices to prove only for $\beta\in]1,+\infty[\setminus\mathbb{N}$ that:\\
\noindent
{\rm(i)} $U_{T_{\beta}}\cap S_c=\emptyset$ for $c<c_{\beta};$\\
\noindent
{\rm(ii)} $U_{T_{\beta}}\cap S_{c_{\beta}}=\lt\{\mu_{\beta}\circ T_{\gamma,1}^{-1}\rt\}_{\gamma\in\mathbb{R}}.$

It is obvious that $c\ge1$ if $U_{T_{\beta}}\cap S_c\neq\emptyset.$

Throughout this proof, all the equations and inequalities for densities hold a.e., and thus we omit \lq\lq a.e.\rq\rq\ for convenience.

Beforehand, we establish the equations for density of a CIUPM for $T_{\beta}$ to be used throughout the proof.

By Proposition~\ref{le3}(i), it suffices to always consider $\mu$ with ${\rm supp}\ \mu\subset[0,c]$ (and thus ${\rm diam}({\rm supp}\ \mu)\le c$), it is easy to verify by Proposition~\ref{le2} and the definition of a CIUPM that $\mu\in U_{T_\beta}$ if and only if its density $\rho_{\mu}$ satisfies the following equations:
\eqb\label{1-1-1}
\sum_{k=0}^{\lf c\rf}\rho_{\mu}(t+k)=1,\ t\in\lt[0,\la c\ra\rt],
\eqe
\eqb\label{1-2-1}
\sum_{k=0}^{\lf c\rf-1}\rho_{\mu}(t+k)=1,\ t\in[\la c\ra,1],
\eqe
        \eqb \label{1-3}\frac{1}{\beta}\sum_{j=0}^{\lf\beta c\rf}\rho_{\mu}\lt(\frac{t+j}{\beta}\rt)=1,\ t\in[0,\la\beta c\ra],\eqe
    \eqb \label{1-4}\frac{1}{\beta}\sum_{j=0}^{\lf\beta c\rf-1}\rho_{\mu}\lt(\frac{t+j}{\beta}\rt)=1,\ t\in[\la\beta c\ra,1]. \eqe
By \eqref{1-1-1} and \eqref{1-2-1}, $0\le\rho_{\mu}(t)\le1$ for $t\in[0,c].$

Note that for $c\le 1+\frac{1}{\beta},$ $$\la\beta c\ra\le\beta+1-\lf\beta c\rf,$$ and thus
\[\frac{t+j}{\beta}\in[c-1,1],\ \forall\ t\in\lt[0,\beta+1-\lf\beta c\rf\rt],\ j=1,\cdots,\lf\beta c\rf-1,\]
from which it follows that
\eqref{1-1-1}-\eqref{1-4} are equivalent to
\eqb \label{1-2}\rho_{\mu}(t)+\rho_{\mu}(t+1)=1,\ t\in[0,c-1],\eqe
\eqb \label{1-1} \rho_{\mu}(t)=1,\ t\in[c-1,1],\eqe
\eqb\label{2-2}    \rho_{\mu}\lt(\frac{t}{\beta}\rt)+\rho_{\mu}\lt(\frac{t+\lf\beta c\rf}{\beta}\rt)=\beta-\lf\beta c\rf+1,\ t\in[0,\la\beta c\ra],\eqe
\eqb\label{2-3}   \rho_{\mu}\lt(\frac{t}{\beta}\rt)=\beta-\lf\beta c\rf+1,\ t\in\lt[\la\beta c\ra,\min\{1,\beta+1-\lf\beta c\rf\}\rt],\eqe
\eqb\label{2-4}    \rho_{\mu}\lt(\frac{t}{\beta}\rt)+\rho_{\mu}\lt(\frac{t+\lf\beta c\rf-1}{\beta}\rt)=\beta-\lf\beta c\rf+2,\ t\in\lt[\min\{1,\beta+1-\lf\beta c\rf\},1\rt].\eqe
By change of variables, \eqref{2-2}-\eqref{2-4} are equivalent to
    \eqb\label{3-2}\rho_{\mu}(t)+\rho_{\mu}\lt(t+\frac{\lf\beta c\rf}{\beta}\rt)=\beta-\lf\beta c\rf+1,\ t\in\lt[0,\frac{\la\beta c\ra}{\beta}\rt],\eqe
   \eqb\label{3-3}\rho_{\mu}(t)=\beta-\lf\beta c\rf+1,\ t\in\lt[\frac{\la\beta c\ra}{\beta},\min\lt\{\frac{1}{\beta},\frac{1}{\beta}+1-\frac{\lf\beta c\rf}{\beta}\rt\}\rt],\eqe
    \eqb\label{3-4}\rho_{\mu}(t)+\rho_{\mu}\lt(t+\frac{\lf\beta c\rf-1}{\beta}\rt)=\beta-\lf\beta c\rf+2,\ t\in\lt[\min\lt\{\frac{1}{\beta},\frac{1}{\beta}+1-\frac{\lf\beta c\rf}{\beta}\rt\},\frac{1}{\beta}\rt].\eqe
In the following, we first prove $U_{T_\beta}\cap S_c=\emptyset$ for $c<c_{\beta}.$

Suppose by way of contradiction that there exists $\mu\in U_{T_\beta}\cap S_c.$ Then its associated density satisfies \eqref{1-2}, \eqref{1-1}, \eqref{3-2}-\eqref{3-4}.
Since $1\le c<c_{\beta}=1+\frac{1}{\beta}$ and $\beta\notin\mathbb{N},$
$$\min\lt\{\frac{1}{\beta},\frac{1}{\beta}+1-\frac{\lf\beta c\rf}{\beta}\rt\}>\frac{\la\beta c\ra}{\beta}.$$
By \eqref{3-3} and $\rho_{\mu}\le1,$we have $\beta\le\lf\beta c\rf,$ yielding that $$\min\lt\{\frac{1}{\beta},\frac{1}{\beta}+1-\frac{\lf\beta c\rf}{\beta}\rt\}=\frac{1}{\beta}+1-\frac{\lf\beta c\rf}{\beta}.$$

Hence \eqref{3-3} and \eqref{3-4} are equivalent to
\[
\rho_{\mu}(t)=\beta-\lf\beta c\rf+1,\ t\in\lt[\frac{\la\beta c\ra}{\beta},\frac{1}{\beta}+1-\frac{\lf\beta c\rf}{\beta}\rt],\]
    \[\rho_{\mu}(t)+\rho_{\mu}\lt(t+\frac{\lf\beta c\rf-1}{\beta}\rt)=\beta-\lf\beta c\rf+2,\ t\in\lt[\frac{1}{\beta}+1-\frac{\lf\beta c\rf}{\beta},\frac{1}{\beta}\rt].\]
Since $\beta\not\in\mathbb{N},$ by $1\le c<1+\frac{1}{\beta}$ and $\beta\le\lf\beta c\rf,$ we have $\lf\beta c\rf=\lf\beta\rf+1.$
This further implies that \eqref{1-2}-\eqref{2-4} are equivalent to
   \eqb\label{7-2} \rho_{\mu}(t)+\rho_{\mu}(t+1)=1,\ t\in[0,c-1],\eqe
    \eqb\label{7-1} \rho_{\mu}(t)=1,\ t\in[c-1,1],\eqe
\eqb\label{7-3} \rho_{\mu}(t)+\rho_{\mu}\lt(t+1+\frac{1-\la\beta\ra}{\beta}\rt)=\la\beta\ra,\ t\in\lt[0,c-1-\frac{1}{\beta}+\frac{\la\beta\ra}{\beta}\rt],\eqe
   \eqb\label{7-4} \rho_{\mu}(t)=\la\beta\ra,\ t\in\lt[c-1-\frac{1}{\beta}+\frac{\la\beta\ra}{\beta},\frac{\la\beta\ra}{\beta}\rt]
   ,\eqe
    \eqb\label{7-5}\rho_{\mu}(t)+\rho_{\mu}\lt(t+1-\frac{\la\beta\ra}{\beta}\rt)=
    \la\beta\ra+1,\ t\in\lt[\frac{\la\beta\ra}{\beta},\frac{1}{\beta}\rt].\eqe
If $c-1<\frac{\la\beta\ra}{\beta},$ \eqref{7-1} contradicts \eqref{7-4} simply because the corresponding intervals have a non-trivial intersection. For the rest of the argument, we assume $c-1\ge\frac{\la\beta\ra}{\beta}.$ Now we aim for a contradiction case by case.\\
  \noindent
    {\rm(i-1)} $\beta\notin\mathbb{Q}.$
Let $L_1:=\lt[c-1-\frac{1}{\beta}+\frac{\la\beta\ra}{\beta},\frac{\la\beta\ra}{\beta}\rt]
.$ Note that\\ $\lambda(L_1)=1+\frac{1}{\beta}-c>0.$ By \eqref{7-2}, \eqref{7-1} and \eqref{7-4}, we have $$\rho_{\mu}(t)=\la\beta\ra,\ t\in L_1,$$$$\rho_{\mu}(t)=1-\la\beta\ra,\ t\in R_1:=L_1+1.$$ Since $$\lt(A_1+\lt(1+\frac{1-\la\beta\ra}{\beta}\rt)\rt)
\bigcup\lt(A_2+\lt(1-\frac{\la\beta\ra}{\beta}\rt)\rt)=[1,c]$$ with $A_1=\lt[0,c-1-\frac{1}{\beta}+\frac{\la\beta\ra}{\beta}\rt]$ and $A_2=\lt[\frac{\la\beta\ra}{\beta},\frac{1}{\beta}\rt],$ by either \eqref{7-3} or \eqref{7-5}, $0\le\rho_{\mu}(t)\le1,$ as well as $\beta\not\in\mathbb{Q},$ we deduce $$\ \ \ \ \ \ \ \ \ \ \ \rho_{\mu}(t)=\lt\la\rho_{\mu}(t)\rt\ra=\lt\la2\lt\la\beta\rt\ra\rt\ra=\lt\la2\beta\rt\ra,\ t\in L_2\subset[0,c-1],$$ where $L_2$ is a union of at most two subintervals of $[0,c-1]$ with $\lambda(L_2)=1+\frac{1}{\beta}-c.$ By induction, we can show that for every $k\in\mathbb{N},$ there exists $L_k,$ a union of finite subintervals  of $[0,c-1]$ with  $\lambda(L_k)=1+\frac{1}{\beta}-c$ such that  \eqb\label{1-5}\rho_{\mu}(t)=\lt\la k\beta\rt\ra,\ t\in L_k.\eqe   Since $\beta\not\in\mathbb{Q},$ $$\lt\la i\beta\rt\ra\neq\lt\la j\beta\rt\ra,\ \forall\ i\neq j,\ i,j\in\mathbb{N},$$ and thus by \eqref{1-5}, $$\lambda(L_i\cap L_j)=0,\ \forall\ i\neq j,\ i,j\in\mathbb{N}.$$ Hence $$\lambda\lt(\cup_{j=1}^kL_j\rt)=k\lt(1+\frac{1}{\beta}-c\rt),\ \forall\ k\in\mathbb{N}.$$ On the other hand, since $\cup_{j=1}^kL_j$ is a subset of $[0,c-1],$ we have $\lambda\lt(\cup_{j=1}^kL_j\rt)\le c-1.$ Take $k=\lt\lfloor\frac{c-1}{1+\frac{1}{\beta}-c}\rt\rfloor+1$ and we arrive at a contradiction.\\
\noindent
{\rm(i-2)} $\beta\in\mathbb{Q}.$ Similarly to case (i-1), for $k=1,\cdots,q_{\beta}-1,$ there exists $L_k\subset[0,c-1]$ with $\lambda(L_k)=1+\frac{1}{\beta}-c$ such that $$\rho_{\mu}(t)=\lt\la k\lt\la\beta\rt\ra\rt\ra=\lt\la\frac{ks_{\beta}}{q_{\beta}}\rt\ra,\ t\in L_k.$$  Since $s_{\beta}$ and $q_{\beta}$ are coprime, by \cite[Chapter~1, Theorem~5.1]{H}, $$\lt\la\frac{is_{\beta}}{q_{\beta}}\rt\ra\neq\lt\la\frac{js_{\beta}}{q_{\beta}}
\rt\ra,\ \forall\ i\neq j,\ 1\le i,j\le q_{\beta}-1,$$ which implies that $$\lambda(L_i\cap L_j)=0,\ \forall\ i\neq j,\ 1\le i,j\le q_{\beta}-1.$$ Hence $\lambda\lt(\cup_{k=1}^{q_{\beta}-1}L_k\rt)=(q_{\beta}-1)\lt(1+\frac{1}{\beta}
-c\rt).$ On the other hand, since $\cup_{k=1}^{q_{\beta}-1}L_k\subset[0,c-1],$  $$(q_{\beta}-1)\lt(1+\frac{1}{\beta}-c\rt)\le c-1,$$ i.e., $c\ge1+\frac{1}{\beta}-\frac{1}{p_{\beta}}=c_{\beta},$ contradicting the assumption that $c<c_{\beta}.$

Next, we show (ii). By definition, it is straightforward to verify that $\mu_{\beta}\in U_{T_{\beta}}\cap S_{c_{\beta}},$ i.e., $\rho_{\lt\la\mu_{\beta}\rt\ra}=\rho_{\lt\la\mu\circ T_{\beta}^{-1}\rt\ra}\equiv1.$

For $\beta\notin\mathbb{Q},$ by \eqref{0-1}, \[\rho_{\lt\la\mu_{\beta}\rt\ra}(t)=\beta t+(-\beta(t+1)+1+\beta)=1,\ t\in\lt[0,1/\beta\rt[;\ \rho_{\la\mu_{\beta}\ra}(t)=1,\ t\in\lt[1/\beta,1\rt[,\]
i.e., $\rho_{\lt\la\mu_{\beta}\rt\ra}\equiv1.$ Note that \[\rho_{\mu_{\beta}\circ T_{\beta}^{-1}}(t)=\begin{cases}
      \frac{1}{\beta}\lt(\beta\frac{t}{\beta}\rt)\ \ \ \ \ \ \ \ \ \ \ \ \ \ \ \ \text{if}\ t\in\lt[0,1\rt[\\
      \frac{1}{\beta}\ \ \ \ \ \ \ \ \ \ \ \ \ \ \ \ \ \ \ \ \ \ \ \ \ \text{if}\ t\in\lt[1,\beta\rt[\\
      \frac{1}{\beta}\lt(-\beta\frac{t}{\beta}+1+\beta\rt)\ \text{if}\ t\in\lt[\beta,\beta+1\rt[\\
      0\ \ \ \ \ \ \ \ \ \ \ \ \ \ \ \ \ \ \ \ \ \ \ \ \ \text{elsewhere}
    \end{cases}=\begin{cases}
      \frac{t}{\beta}\ \ \ \ \ \ \ \ \ \ \ \ \ \ \ \ \ \text{if}\ t\in\lt[0,1\rt[,\\
      \frac{1}{\beta}\ \ \ \ \ \ \ \ \ \ \ \ \ \ \ \ \ \text{if}\ t\in\lt[1,\beta\rt[,\\
      -\frac{t}{\beta}+1+\frac{1}{\beta}\ \text{if}\ t\in\lt[\beta,\beta+1\rt[,\\
      0\ \ \ \ \ \ \ \ \ \ \ \ \ \ \ \ \  \text{elsewhere},
    \end{cases}
\] and thus \[\rho_{\lt\la\mu_{\beta}\circ T_{\beta}^{-1}\rt\ra}(t)=\frac{t}{\beta}+\lf\beta\rf\cdot\frac{1}{\beta}+1+\frac
{1}{\beta}-\frac{t+\lf\beta\rf+1}{\beta}=1,\ t\in\lt[0,\lt\la\beta\rt\ra\rt[;\] \[\rho_{\lt\la\mu_{\beta}\circ T_{\beta}^{-1}\rt\ra}(t)=\frac{t}{\beta}+\lt(\lf\beta\rf-1\rt)\cdot\frac{1}
{\beta}+1+\frac{1}{\beta}-\frac{t+\lf\beta\rf}{\beta}=1,\ t\in\lt[\la\beta\ra,1\rt[,\]
i.e., $\rho_{\lt\la\mu_{\beta}\circ T_{\beta}^{-1}\rt\ra}\equiv1.$ Thus $\mu_{\beta}\in U_{T_{\beta}}\cap S_{c_{\beta}}.$

For $\beta\in\mathbb{Q},$ by \eqref{0-3} and induction, it is easy to confirm that \[\begin{split}&\rho_{\lt\la\mu_{\beta}\rt\ra}(t)=\frac{j}{q_{\beta}}+\frac{q_
{\beta}-j}{q_{\beta}}=1,\ t\in\lt[\frac{j-1}{p_{\beta}},\frac{j}{p_{\beta}}\rt[,\ j=1,\cdots,q_{\beta}-1;\\ &\rho_{\lt\la\mu_{\beta}\rt\ra}(t)=1,\ t\in\lt[\frac{q_{\beta}-1}{p_{\beta}},1\rt[,\end{split}\]
i.e., $\rho_{\lt\la\mu_{\beta}\rt\ra}\equiv1.$ Again by induction, one can show for $j=1,\cdots,s_{\beta},$
\[\rho_{\mu_{\beta}\circ T_{\beta}^{-1}}(t)=\begin{cases}
      \frac{j}{p_{\beta}}\ \ \ \text{if}\ t\in\lt[\frac{j-1}{q_{\beta}},\frac{j}{q_{\beta}}\rt[,\\
      \frac{q_{\beta}}{p_{\beta}}\ \ \ \text{if}\ t\in\lt[k+\frac{j-1}{q_{\beta}},k+\frac{j}{q_{\beta}}\rt[,\ \text{for}\ k=1,\cdots,\lf\beta\rf,\\
      \frac{s_{\beta}-j}{p_{\beta}}\ \text{if}\ t\in\lt[\lf\beta\rf+1+\frac{j-1}{q_{\beta}},\lf\beta\rf+1+\frac{j}{q_
      {\beta}}\rt[,
    \end{cases}\] and for $j=s_{\beta}+1,\cdots,q_{\beta},$\[\rho_{\mu_{\beta}\circ T_{\beta}^{-1}}(t)=\begin{cases}
      \frac{j}{p_{\beta}}\ \ \ \ \ \ \ \  \text{if}\ t\in\lt[\frac{j-1}{q_{\beta}},\frac{j}{q_{\beta}}\rt[,\\
      \frac{q_{\beta}}{p_{\beta}}\ \ \ \ \ \ \ \ \text{if}\ t\in\lt[k+\frac{j-1}{q_{\beta}},k+\frac{j}{q_{\beta}}\rt[,\ \text{for}\ k=1,\cdots,\lf\beta\rf-1,\\
      \frac{q_{\beta}+s_{\beta}-j}{p_{\beta}}\ \text{if}\ t\in\lt[\lf\beta\rf+\frac{j-1}{q_{\beta}},\lf\beta\rf+\frac{j}{q_
      {\beta}}\rt[,
    \end{cases}\] yielding
    \[\rho_{\lt\la\mu_{\beta}\circ T_{\beta}^{-1}\rt\ra}(t)=\frac{j}{p_{\beta}}+\lf\beta\rf\cdot\frac{q_
    {\beta}}{p_{\beta}}+\frac{s_{\beta}-j}{p_{\beta}}=1,\ t\in\lt[\frac{j-1}{q_{\beta}},\frac{j}{q_{\beta}}\rt[,\ \text{for}\ j=1,\cdots,s_{\beta},\]\[\rho_{\lt\la\mu_{\beta}\circ T_{\beta}^{-1}\rt\ra}(t)=\frac{j}{p_{\beta}}+(\lf\beta\rf-1)\cdot\frac{
    q_{\beta}}{p_{\beta}}+\frac{q_{\beta}+s_{\beta}-j}{p_{\beta}}=1,\ t\in\lt[\frac{j-1}{q_{\beta}},\frac{j}{q_{\beta}}\rt[,\ \text{for}\ j=s_{\beta}+1,\cdots,q_{\beta},\]i.e., $\rho_{\lt\la\mu_{\beta}\circ T_{\beta}^{-1}\rt\ra}\equiv1.$

Thus, by Proposition~\ref{le3} (i), it suffices to show that
\begin{claim}
If $\mu\in U_{T_{\beta}}\cap S_{c_{\beta}}$ with ${\rm supp}\ \mu\subset[0,c_{\beta}],$ then $\mu=\mu_{\beta}.$
\end{claim}
In the following, we prove this claim case by case.\\
\noindent
{\rm(ii-1)} $\beta\notin\mathbb{Q}.$
In this case, it seems not enough to only deal with equations and inequalities for the density (which only holds in the almost everywhere sense); we instead need to consider the distribution function. Recall that $F_{\mu}$ is continuous for all $\mu\in U_{T_{\beta}},$ by Proposition~\ref{le4}.

It follows from \eqref{1-2}, \eqref{1-1}, \eqref{3-2}-\eqref{3-4} together with the continuity of $F_{\mu}$ that, for $c=1+\frac{1}{\beta},$
    \eqb\label{8-2} F_{\mu}(t)+F_{\mu}(t+1)=t+F_{\mu}(1),\ t\in\lt[0,\frac{1}{\beta}\rt],\eqe
   \eqb\label{8-1} F_{\mu}(t)=F_{\mu}\lt(\frac{1}{\beta}\rt)+t-\frac{1}{\beta},\ t\in\lt[\frac{1}{\beta},1\rt],\eqe
 \eqb\label{8-3} F_{\mu}(t)+F_{\mu}\lt(t+1+\frac{1-\la\beta\ra}{\beta}\rt)=\la\beta\ra t+F_{\mu}
 \lt(1+\frac{1-\la\beta\ra}{\beta}\rt),\ t\in\lt[0,\frac{\la\beta\ra}{\beta}\rt],\eqe
    \eqb\label{8-4}\begin{split} &F_{\mu}(t)+F_{\mu}\lt(t+1-\frac{\la\beta\ra}{\beta}\rt)\\=&(\la\beta\ra+1)
    \lt(t-\frac{\la\beta\ra}{\beta}\rt)+F_{\mu}\lt(\frac{\la\beta\ra}{\beta}\rt)+
    F_{\mu}(1),\ t\in\lt[\frac{\la\beta\ra}{\beta},\frac{1}{\beta}\rt].\end{split}\eqe
By \eqref{8-2} and \eqref{8-3},
\eqb\label{9-1}
\frac{F_{\mu}\lt(t+\frac{1-\la\beta\ra}{\beta}\rt)-F_{\mu}(t)}{\frac{1-\la\beta\ra}
{\beta}}-\beta t=C_1,\ t\in\lt[0,\frac{\la\beta\ra}{\beta}\rt]
\eqe    with $C_1=\frac{F_{\mu}\lt(\frac{1}{\beta}\rt)-F_{\mu}\lt(\frac{\la\beta\ra}{\beta}
\rt)}{\frac{1-\la\beta\ra}{\beta}}-\la\beta\ra.$ Similarly,
by \eqref{8-2} and \eqref{8-4},
\eqb\label{9-2}
\frac{F_{\mu}\lt(t+\frac{\la\beta\ra}{\beta}\rt)-F_{\mu}(t)}{\frac{\la\beta\ra}
{\beta}}-\beta t=C_2,\ t\in\lt[0,\frac{1-\la\beta\ra}{\beta}\rt]
\eqe    with $C_2=\frac{F_{\mu}\lt(\frac{\la\beta\ra}{\beta}\rt)}{\frac{\la\beta\ra}{\beta}}.$

Furthermore,  by \eqref{9-1} and \eqref{9-2}, we can show by induction that for all $m,n\in\mathbb{Z}$ satisfying $m\frac{\la\beta\ra}{\beta}+n\frac{1-\la\beta\ra}{\beta}\in\lt]0,\frac{1}{\beta}
\rt[,$\eqb
\label{5-1-1}
\begin{split}
&F_{\mu}\lt(m\frac{\la\beta\ra}{\beta}+n\frac{1-\la\beta\ra}{\beta}\rt)\\=&\frac
{\beta}{2}\lt(m\frac{\la\beta\ra}{\beta}+n\frac{1-\la\beta\ra}{\beta}\rt)^2
+\lt(C_1-\frac{1-\la\beta\ra}{2}\rt)n\frac{1-\la\beta\ra}{\beta}\\&+\lt(C_2-\frac
{\la\beta\ra}{2}\rt)m\frac{\la\beta\ra}{\beta}.
\end{split}\eqe
By Proposition~\ref{le1}, $\lt(m\frac{\la\beta\ra}{\beta}+n\frac{1-\la\beta\ra}{\beta}\rt)_{m,n\in\mathbb{Z}}
\bigcap\lt]0,\frac{1}{\beta}\rt[$ is dense in $\lt[0,\frac{1}{\beta}\rt].$ Thus, for every $t\in\lt[0,\frac{1}{\beta}\rt]\lt\backslash\lt(m\frac{\la\beta\ra}{\beta}+n\frac
{1-\la\beta\ra}{\beta}\rt)_{m,n\in\mathbb{Z}},\rt.$ there exist two sequences $(m_k)_{k\in\mathbb{N}}$ and $(n_k)_{k\in\mathbb{N}}$ such that $$\lim_{k\to\infty}m_k\frac{\la\beta\ra}{\beta}+n_k\frac{1-\la\beta\ra}{\beta}=t.$$ It is easy to see that $\lim_{k\to\infty}|m_k|=\lim_{k\to\infty}|n_k|=+\infty$ (otherwise, both $(m_k)_{k\in\mathbb{N}}$ and $(n_k)_{k\in\mathbb{N}}$ are bounded, and thus $t\in\lt(m\frac{\la\beta\ra}{\beta}+n\frac{1-\la\beta\ra}{\beta}\rt)_{m,n\in
\mathbb{Z}}$\Big). Substituting $(m,n)$ in \eqref{5-1-1} by $(m_k,n_k)$ and letting $k\to\infty$ on both sides of \eqref{5-1-1}, by the continuity of $F_{\mu},$ $$ C_1-\frac{1-\la\beta\ra}{2}=C_2-\frac{\la\beta\ra}{2}.$$ From \eqref{5-1-1} it follows that \eqb\label{5-2}
F_{\mu}(t)=\frac{\beta}{2}t^2+Ct,\ \forall\ t\in\lt[0,\frac{1}{\beta}\rt],
\eqe where $C=C_1-\frac{1-\la\beta\ra}{2}.$ By the definition of derivative, it follows from \eqref{5-2} that
$$F_{\mu}'(t)=\beta t+C, \forall\ t\in\lt]0,\frac{1}{\beta}\rt[.$$ Since $F_{\mu}$ is non-decreasing in $\lt]0,\frac{1}{\beta}\rt[,$ $\lim_{t\downarrow0}F_{\mu}'(t)\ge0$ implies that $C\ge0.$ By \eqref{8-2} and \eqref{5-2}, \[
F_{\mu}(t)=-\frac{\beta}{2}(t-1)^2+(1-C)(t-1)+F_{\mu}(1),\ \forall\ t\in\lt[1,1+\frac{1}{\beta}\rt].
\] Similarly, $\lim_{t\uparrow\left(1+\frac{1}{\beta}\right)}F_{\mu}'(t)\ge0$ yields $C\le0.$ Thus $C=0.$

By \eqref{8-1}, $F_{\mu}$ is given by
\[ F_{\mu}(t)=\cab
\frac{\beta}{2}t^2\ \ \ \ \ \ \ \ \ \ \ \ \ \ \ \ \ \ \ \ \ \ \ \ \ \ \text{if}\ t\in\lt[0,\frac{1}{\beta}\rt[,\\ t-\frac{1}{2\beta}\ \ \ \ \ \ \ \ \ \ \ \ \ \ \ \ \ \ \ \ \ \ \text{if}\ t\in\lt[\frac{1}{\beta},1\rt[,\\-\frac{\beta}{2}(t-1)^2+t-\frac{1}{2\beta}\ \text{if}\ t\in\lt[1,1+\frac{1}{\beta}\rt[,
\cae
\] equivalently, $\rho_{\mu}=\rho_{\beta}$ and thus $\mu=\mu_{\beta}.$\\
\noindent
{\rm(ii-2)} $\beta\in\mathbb{Q}.$ Recall the definitions of $p_{\beta},\ q_{\beta}$ and $s_{\beta}$ for every $\beta\in\mathbb{Q}\setminus\{0\}$ in the previous section, we know $\lf\beta\rf=\lt\lf\frac{p_{\beta}-1}{q_{\beta}}\rt\rf$ for $\beta\notin\mathbb{N}.$ Hence
\eqb\label{11-2}
\rho_{\mu}(t)+\rho_{\mu}(t+1)=1, t\in\lt[0,\frac{q_{\beta}-1}{p_{\beta}}\rt],
\eqe\eqb\label{11-1}
\rho_{\mu}(t)=1,\ t\in\lt[\frac{q_{\beta}-1}{p_{\beta}},1\rt],
\eqe
\eqb\label{11-3}
\rho_{\mu}(t)+\rho_{\mu}\lt(t+1+\frac{q_{\beta}-s_{\beta}}{p_{\beta}}\rt)
=\frac{s_{\beta}}{q_{\beta}},\ t\in\lt[0,\frac{s_{\beta}-1}{p_{\beta}}\rt],
\eqe
\eqb\label{11-4}
\rho_{\mu}(t)=\frac{s_{\beta}}{q_{\beta}},\ t\in\lt[\frac{s_{\beta}-1}{p_{\beta}},\frac{s_{\beta}}{p_{\beta}}\rt],
\eqe
\eqb\label{11-5}
\rho_{\mu}(t)+\rho_{\mu}\lt(t+1-\frac{s_{\beta}}{p_{\beta}}\rt)=1+
\frac{s_{\beta}}{q_{\beta}},\ t\in\lt[\frac{s_{\beta}}{p_{\beta}},\frac{q_{\beta}}{p_{\beta}}\rt].
\eqe
It follows from \eqref{11-3} and \eqref{11-5} that, \[\lt(\lt[0,\frac{s_{\beta}-1}{p_{\beta}}\rt[+\lt(1+\frac{q_{\beta}-s_{\beta}}
{p_{\beta}}\rt)\rt)\bigcup\lt(\lt
[\frac{s_{\beta}}{p_{\beta}},\frac{q_{\beta}}{p_{\beta}}\rt[+\lt(1-\frac{s_
{\beta}}{p_{\beta}}\rt)\rt)=\lt[1,1+\frac{q_{\beta}-1}{p_{\beta}}\rt[.\]
Using \eqref{11-2}, \eqref{11-3} and \eqref{11-5},
\[\rho_{\mu}\lt(t+\frac{s_{\beta}}{p_{\beta}}\rt)-\rho_{\mu}(t)=\frac{s_
{\beta}}{q_{\beta}},\ t\in\lt[0,\frac{q_{\beta}-s_{\beta}}{q_{\beta}}\rt],\]
\[\rho_{\mu}\lt(t+\frac{q_{\beta}-s_{\beta}}{p_{\beta}}\rt)-\rho_{\mu}(t)=1-
\frac{s_{\beta}}{q_{\beta}},\ t\in\lt[0,\frac{s_{\beta}-1}{q_{\beta}}\rt].\]

Similarly to \eqref{5-1-1}, we can show by induction that
\eqb\label{13-1}
\rho_{\mu}\lt(t+m\frac{s_{\beta}}{p_{\beta}}+n\frac{q_{\beta}-s_{\beta}}{p_
{\beta}}\rt)=\rho_{\mu}(t)+m\frac{s_{\beta}}{q_{\beta}}+n\lt
(1-\frac{s_{\beta}}{q_{\beta}}\rt),
\eqe for $m,n\in\mathbb{Z},$ $t\in\lt[0,\frac{q_{\beta}-1}{p_{\beta}}\rt]$ a.e. satisfying $t+m\frac{s_{\beta}}{p_{\beta}}+n\frac{q_{\beta}-s_{\beta}}{p_{\beta}}\in
\lt[0,\frac{q_{\beta}-1}{p_{\beta}}\rt].$
Since $s_{\beta}$ and $q_{\beta}$ are coprime, from for instance \cite[Chapter~1, Theorem~4.4(i)]{H}, there exist $m_0,n_0\in\mathbb{Z}$ such that $m_0s_{\beta}+n_0(q_{\beta}-s_{\beta})=1.$  Then it follows from \eqref{13-1} that
\eqb
\label{13-2}
\rho_{\mu}\lt(t+\frac{j}{p_{\beta}}\rt)=\rho_{\mu}(t)+\frac{j}{q_{\beta}},
\eqe for $j\in\mathbb{Z},$ $t\in\lt[0,\frac{q_{\beta}-1}{p_{\beta}}\rt]$ a.e. satisfying $t+\frac{j}{p_{\beta}}\in\lt[0,\frac{q_{\beta}-1}{p_{\beta}}\rt].$  By \eqref{11-1}, \eqref{11-2}, \eqref{11-4} and \eqref{13-2}, we can prove by induction that $\rho_{\mu}=\rho_{\beta}$ and thus $\mu=\mu_{\beta}.$
\pre
\rb
{\rm (i)} For $\beta\neq0,$ $\alpha\in\mathbb{R},$ it follows from Theorem~\ref{th1} that there always exists a CIUPM for $T_{\alpha,\beta}$ with arbitrary length  (in diameter) $c\ge c_{\beta}.$
Moreover, from the proof of Theorem~\ref{th1} one easily observes that if $\beta\in\mathbb{Q}\cap\ [1,+\infty[,$ then $\widetilde{\mu}_{\beta}$ with its density function
\[\rho_{\widetilde{\mu}_\beta}(t)=\begin{cases}
      \beta t\ \ \ \ \ \ \ \ \ \ \ \ \ \ \ \ \text{if}\ t\in\lt[0,\frac{1}{\beta}\rt[,\\
      1,\ \ \ \ \ \ \ \ \ \ \ \ \ \ \ \ \text{if}\ t\in\lt[\frac{1}{\beta},1\rt[,\\
      -\beta t+1+\beta\ \text{if}\ t\in\lt[1,1+\frac{1}{\beta}\rt],\\
      0\ \ \ \ \ \ \ \ \ \ \ \ \ \ \ \ \ \ \text{elsewhere},
    \end{cases}\] is another CIUPM but with ${\rm diam} \lt(\rm supp\ \widetilde{\mu}_{\beta}\rt)>c_{\beta}.$ \\
\noindent
{\rm (ii)} From the proof of Corollary~\ref{co1}, for every $n\in\mathbb{N}\setminus\{1\},$ every probability vector $(p_1,\cdots,p_n)$ with $p_i>0$ for all $i,$ and $0=q_1<q_2<\cdots<q_n=1,$ for all $c>c_{\beta},$ $\mu=\sum_{i=1}^np_i\mu_{\beta}\circ T_{{q_i(c-c_{\beta}),1}}^{-1}\in U_{T_{\alpha,\beta}}\cap\ S_c.$ This illustrates the non-uniqueness of CIUPMs and reflects the potential complexity of the set $U_{T_{\alpha,\beta}}.$ It may be interesting to completely characterize $U_{T_{\alpha,\beta}}.$
\re
\rb
Notice that $U_{T_{\alpha,\beta}}\cap S_{c_{\beta}}$ may not contain all \lq\lq slimmest\rq\rq\ CIUPMs in the sense that the support has the smallest {\em Lebesgue measure} (instead of the diameter). In other words, there may exist CIUPMs with the smallest {\em disconnected} support in diameter:    $U_{T_{\alpha,\beta}}\cap S_{c_{\beta}}\subsetneqq U_{T_{\alpha,\beta}}\cap \overline{S}_{c_{\beta}}$ for some $\beta\neq0,$ where $\overline{S}_c:=\lt\{\mu\in\mathcal{P}(\mathbb{R}):\ \lambda({\rm supp}(\mu))\le c\rt\}.$ For instance, $\mu:=\lambda|_{[0,1/2]}+\lambda|_{[3/2,2]}$ with $\lambda({\rm supp}\ \mu)=1=c_k$ is a CIUPM for every linear transformation $T_{\alpha,k}$ with $\alpha\in\mathbb{R}$ and nonzero integer $k.$  In fact,
for every pair $(m,n)\in\mathbb{Z}^2,$ every $\beta\in[1,+\infty[\ \cap\ \mathbb{Q},$ define \[\rho_{\beta,m,n}(t)=\begin{cases}
      \frac{j}{q_{\beta}}\ \text{if}\ t\in\lt[mq_{\beta}+\frac{j-1}{p_{\beta}},mq_{\beta}+\frac{j}{p_
      {\beta}}\rt[\bigcup\lt[mq_{\beta}+1+\frac{q_{\beta}-j-1}{p_{\beta}}
      ,mq_{\beta}+1+\frac{q_{\beta}-j}{p_{\beta}}\rt[,\\ \ \ \ \ \ \ \text{for}\ j=1,\cdots,q_{\beta}-1,\\
      1\ \ \text{if}\ t\in\lt[nq_{\beta}+\frac{q_{\beta}-1}{p_{\beta}},nq_{\beta}+1\rt[,\\
      0\ \ \text{elsewhere}.
    \end{cases}\]
It can be shown by induction (analogous to the proof of Theorem~\ref{th1}) that $\mu_{\beta,m,n},$ with density function $\rho_{\beta,m,n},$ is a CIUPM for $T_{\alpha,\beta}$ with $\lambda\lt({\rm supp}\ \mu_{\beta,m,n}\rt)=c_{\beta};$ moreover, $\mu_{\beta,m,n}\in S_{c_{\beta}}$ if and only if $m=n$ when ${\rm supp}\ \mu_{\beta,m,n}$ is an interval. Thus $U_{T_{\alpha,\beta}}\cap S_{c_{\beta}}\subsetneqq U_{T_{\alpha,\beta}}\cap \overline{S}_{c_{\beta}}$ for every $\beta\in\mathbb{Q}\setminus\{0\},$ by Proposition~\ref{le3}. However, due to the nature of irrationality, the author conjectures $U_{T_{\alpha,\beta}}\cap S_{c_{\beta}}=U_{T_{\alpha,\beta}}\cap \overline{S}_{c_{\beta}}$ for every $\beta\in\mathbb{R}\setminus\mathbb{Q}.$
\re

\section*{Acknowledgement}
Deepest thanks to my thesis advisor Arno Berger for proposing this problem, and many helpful discussions and constant encouragement. The author is also indebted to an anonymous referee for his/her proofreading as well as valuable suggestions which help improve the presentation of the manuscript. This research is supported in part by a Pacific Institute for the Mathematical Sciences (PIMS) Graduate Scholarship and a Josephine Mitchell Graduate Scholarship. 

\end{document}